\documentclass{amsart}

\usepackage{graphicx}

\newtheorem{theorem}{Theorem}[section]
\newtheorem{lemma}[theorem]{Lemma}

\theoremstyle{definition}

\theoremstyle{remark}

\numberwithin{equation}{section}

\begin{document}

\title{Global in time existence of strong solution to 3D Navier-Stokes equations}
\author{Abdelkerim Chaabani}
\address{Department of mathematics, Faculty of Mathematical, Physical and Natural \\ Sciences of Tunis, University of Tunis El Manar, Tunis 2092, Tunisia.}
\email{abdelkerim.chaabani@fst.utm.tn}



\begin{abstract}
The purpose of this paper is to bring to light a method through which the global in time existence for arbitrary large in $H^1$ initial data of a strong solution to 3D periodic Navier-Stokes equations follows. The method consists of subdividing the time interval of existence into smaller sub-intervals carefully chosen. These sub-intervals are chosen based on the hypothesis that for any wavenumber $m$, one can find an interval of time on which the energy quantized in low-frequency components (up to $m$) of the solution u is lesser than the energy quantized in high-frequency components (down to $m$) or otherwise the opposite. We associate then a suitable number $m$ to each one of the intervals and we prove that the norm $\|u(t)\|_{\dot{H}^1}$ is bounded in both mentioned cases. The process can be continued until reaching the maximal time of existence $T_{max}$ which yields the global in time existence of strong solution.
\end{abstract}


\maketitle
\section{Introduction}
Let us consider the following incompressible Navier-Stokes equations:
$$
\left\{
  \begin{array}{l}
     \partial_{t}u-\nu\Delta u+(u\cdot\nabla)u+\nabla p=0, \ \ \ (x,t)\in \mathbb{T}^3\times\mathbb{R}_+\\
     \nabla\cdot u=0, \ \ \ (x,t)\in \mathbb{T}^3\times\mathbb{R}_+\\
     u|_{t=0}=u_0(x), \ \ \ x\in\mathbb{T}^3,\\
  \end{array}
\right.
\leqno{(NSE)}
$$
where the constant $\nu>0$ is the viscosity of the fluid, and  $\mathbb{T}^3=\mathbb{R}^3/\mathbb{Z}^3$ is the three-dimensional torus with periodic boundary conditions.
Here $u$ is a three-dimensional vector field $u = (u_1, u_2, u_3)$ representing the velocity of the fluid, and $p$ is a scalar denoting the pressure, both are unknown functions
of the space variable $x$ and time variable $t$. We recall that the pressure can be
eliminated by projecting $(NSE)$ onto the space of free divergence vector fields, using the Leray projector $$\mathbb{P} =Id- \nabla\Delta^{-1}\nabla\cdot.$$
Thus, it will be convenient using the following equivalent system
$$
\left\{
  \begin{array}{l}
     \partial_{t}u-\nu\Delta u+\mathbb{P}(u\cdot\nabla)u=0, \ \ \ (x,t)\in \mathbb{T}^3\times\mathbb{R}_+\\
     \nabla\cdot u=0, \ \ \ (x,t)\in \mathbb{T}^3\times\mathbb{R}_+\\
     u|_{t=0}=u_0(x), \ \ \ x\in\Bbb{T}^3.\\
  \end{array}
\right.
\leqno{(NS)}
$$
We define the Sobolev spaces $H^s(\mathbb{T}^3)$ for $s\geq 0$ by the Fourier expansion $$H^s(\mathbb{T}^3)=\left\{ u\in L^2(\mathbb{T}^3): ~u(x):=\sum_{k\in\mathbb{Z}^3}\hat{u}(k,t)e^{ikx},~~\hat{u}(k,t)=\overline{\hat{u}(k,t)},~\|u\|_{H^s}<\infty \right\},$$
where $$\|u\|_{H^s}^2:=\sum_{k\in\mathbb{Z}^3}(1+|k|^{2s})|\hat{u}(k,t)|^2$$
and $$\hat{u}(k,t)=\int_{\mathbb{T}^3}u(x)e^{-ikx}dx.$$
We also give the definition of homogeneous Sobolev space: $$\dot{H}^s(\mathbb{T}^3)=\left\{ u\in L^2(\mathbb{T}^3): ~u(x):=\sum_{k\in\mathbb{Z}^3}\hat{u}(k,t)e^{ikx},~~\hat{u}(k,t)=\overline{\hat{u}(k,t)},~\|u\|_{\dot{H}^s}<\infty \right\},$$ and endowed by the norm $$\|u\|_{\dot{H}^s}:=\|\Lambda^su\|_{L^2}=\left(\sum_{k\in \mathbb{Z}^3}|k|^{2s}|\hat{u}(k,t)|^2\right)^{1/2}, $$
where by $\Lambda$ we refer to the operator $\sqrt{-\Delta}$.\\
We will also use the following function spaces:
$$\mathcal{D}_\sigma:=\{\varphi\in[C^\infty_c(\mathbb{T}^3)]^3:~\nabla\cdot \varphi=0\}$$
$$L^2_{\sigma}(\mathbb{T}^3):=\mbox{closure of}~\mathcal{D}_\sigma~\mbox{in }L^2$$
$$H^1_{\sigma}(\mathbb{T}^3):=\mbox{closure of}~\mathcal{D}_\sigma~\mbox{in }H^1$$

For an initial data $u_0\in L^2_{\sigma}(\mathbb{T}^3)$, it was proven by Leray and Hopf  that there exists a global weak solution $u\in L^\infty_t(L^2_\sigma)\cap L^2_t(H^1_\sigma)$.
\begin{theorem}
 For every $u_0\in L^2_{\sigma}(\mathbb{T}^3)$ there exists at least one global in time weak solution $u\in L^\infty(0,\infty;L^2_{\sigma}(\mathbb{T}^3))\cap L^2(0,\infty;H^1_{\sigma}(\mathbb{T}^3)) $ of the Navier-Stokes equations satisfying the initial condition $u_0$. In particular, $u$ satisfies the energy inequality
 \begin{eqnarray}\label{W.E}
 \frac{1}{2}\|u(t)\|_{L^2_{\sigma}(\mathbb{T}^3)}^2+\nu\int_0^t\|\nabla u(\tau)\|_{L^2_{\sigma}(\mathbb{T}^3)}^2d\tau\leq \frac{1}{2}\|u_0\|_{L^2_{\sigma}(\mathbb{T}^3)}^2.
 \end{eqnarray}
\end{theorem}
This result was proved by \cite{Hopf} as a generalisation of a previous existence theorem due to  \cite{Leray} for the whole space $\mathbb{R}^3$.\\
It is also known that local in time strong solutions exist on the whole space due to  \cite{Leray}, while the case of a bounded domain is due to \cite{kis}.
\begin{theorem}
 There is a constant $C>0$ such that any initial condition $u_0\in H_{\sigma}^1(\mathbb{T}^3)$ gives rise to a strong solution of the Navier-Stokes equations $$u\in L^\infty(0,T_{max};H^1_\sigma)\cap L^2(0,T_{max};H^2),~~\mbox{where~} T_{max}= \frac{C}{\|\nabla u_0\|_{L^2}^4}.$$
\end{theorem}
The existence of global in time strong solution is known to occur for small initial data due to \cite{kato} and \cite{chemin}. However, it remains the major open problem as to whether these solutions can be extended to be global in time for arbitrary large in $H^1$ initial data. Originally, the problem is the question of global existence of smooth solutions to the Navier-Stokes equations satisfying bounded energy condition (i.e.: $u\in C^\infty( \mathbb{T}^3\times\mathbb{R}_+)$ and $\int_{\mathbb{T}^3}|u(x)|^2dx<\infty$)  or otherwise a breakdown. The official description has been given by Fefferman in \cite{Feff}. The official Clay Millennium problem is to give a proof of one of the four following statements:
\begin{itemize}
    \item[(A)]  Existence and smoothness of Navier–Stokes solutions on $\mathbb{R}^3$
    \item[(B)]  Existence and smoothness of Navier–Stokes solutions on $\mathbb{R}^3/\mathbb{Z}^3$
    \item[(C)] Breakdown of Navier–Stokes solutions on $\mathbb{R}^3$
    \item[(C)] Breakdown of Navier–Stokes solutions on $\mathbb{R}^3/\mathbb{Z}^3$
\end{itemize}
In this paper, we prove the statement (B) which can be alternatively formulated as follows:
\begin{theorem}\label{main}
For every $u_0\in H^1_{\sigma}(\mathbb{T}^3)$ there exists a unique  global in time strong solution $u\in L^\infty(0,\infty;H^1_{\sigma}(\mathbb{T}^3))\cap L^2(0,\infty;H^2(\mathbb{T}^3))$ of the Navier-Stokes equations.
\end{theorem}
The method used to extend the solution into a global one is to prove that on an interval of strictly positive length $[t_0,t_1]\subset (0,T_{max})$ under a first condition on $\sum_{k\in\mathbb{Z}^3}|\hat{u}(k,t)|$ among two possible ones: $$\underbrace{\sum_{|k|> m}|\hat{u}(k,t)|\leq \sum_{|k|\leq  m}|\hat{u}(k,t)|}_{\mbox{condition 1}}~~\mbox{or~}\underbrace{\sum_{|k|\leq m}|\hat{u}(k,t)|\leq \sum_{|k|> m}|\hat{u}(k,t)|}_{\mbox{condition 2}},$$ the solution will be controlled in $H^1$ by a suitable function defined in terms of time $t$, $\|\nabla u(t_0)\|_{L^2}$, $\| u_0\|_{L^2}$ and a finite number $m$ (depending on the viscosity $\nu$ and $\|\nabla u(t_0)\|_{L^2}$), until reaching $t_1$.
Otherwise, that is if condition 2 holds true $\forall t\in [t_0,t_1]$, then the norm $\|\nabla u(t)\|_{L^{2}}^2$ is non-increasing on $[t_0,t_1]$. We continue then in this vein until reaching $T_{max}.$ To be more precise, we subdivide the interval $(0,T_{max})$ into a series of successive sub-intervals each of them is akin to $[t_0,t_1]$, i.e.: on each of them either condition $1$ or $2$ holds. It should be emphasized that the number $m$ may change from an interval to another.  \\
We quote the following two results, the proof of which is given in \cite{james} and based on that of (Theorem 10.6 \cite{Cons}).
\begin{theorem}\label{regularity}
Let $u$ be a strong solution of the Navier-Stokes equations $(NS)$ on the time interval $[0, T]$, with initial condition $u_0\in H^1$. Then for all $0<\varepsilon<T$ we have
$u \in C([\varepsilon, T];H^p) \mbox{~~for all~} p \in \mathbb{N}.$
In particular, for all $t\in [0,T]$ the function $u(t)$ is smooth with respect to the space variables.
\end{theorem}
The Theorem \ref{regularity} together with the following lemma constitute a cornerstone in establishing the proof of Theorem \ref{main}.
\begin{lemma}\label{regularity2}
Let $u$ be a strong solution of the Navier-Stokes equations on the
time interval $[0,T]$. Then for every $\varepsilon>0$, and all $p,~l\in \mathbb{N}$ we have $$\partial^l_tu\in L^\infty(\varepsilon;T;H^p).$$
\end{lemma} 
The rest of the paper is dedicated to give the proof of Theorem \ref{main}.
\section{The proof}
The analysis can be started by sketching the procedure through which the existence of local in time strong solution to $(NS)$ follows. To this end, let $P_n$ be the projection onto the Fourier modes of order up to $n\in \mathbb{N}$, that is $$P_n(\sum_{k\in\mathbb{Z}^3}\hat{\vartheta}_ke^{ixk})=\sum_{|k|\leq n}\hat{\vartheta}_ke^{ixk}.$$
Let $u_n=P_nu$ be the solution to
$$
\left\{
  \begin{array}{l}
     \partial_{t}u_n-\nu\Delta u_n+P_n[(u_n\cdot\nabla)u_n]=0, \ \ \ (x,t)\in \mathbb{T}^3\times\mathbb{R}_+\\
     \nabla\cdot u_n=0, \ \ \ (x,t)\in \mathbb{T}^3\times\mathbb{R}_+\\
     u_n|_{t=0}(x)=P_n(u^{in})(x), \ \ \ x\in\Bbb{T}^3.\\
  \end{array}
\right.
\leqno{(NS_n)}
$$
For some $T_n$, there exists a  solution $u_n\in C^\infty([0,T_n)\times\mathbb{T}^3)$ to this finite-dimensional locally-Lipschitz system of ODEs. We take the $L^2$-inner product of the first equation in ($NS_n$) against $-\Delta u_n$ to obtain
\begin{eqnarray}
\frac{1}{2}\frac{d}{dt}\| \nabla u_n(t)\|_{L^{2}}^2+\nu\| \Delta u_n(t)\|_{L^{2}}^2&\leq&
 |\langle (u_n\cdot\nabla u_n), \Delta u_n \rangle_{L^2(\mathbb{T}^3)}|\nonumber\\
&\leq& \|u_n(t)\|_{L^\infty(\mathbb{T}^3)}\|\nabla u_n(t)\|_{L^2(\mathbb{T}^3)}\|\Delta u_n(t)\|_{L^2(\mathbb{T}^3)}\nonumber\\
&\leq& c\| u_n(t)\|_{H^1(\mathbb{T}^3)}^{1/2}\| u_n(t)\|_{H^2(\mathbb{T}^3)}^{1/2}\|\nabla u_n(t)\|_{L^2(\mathbb{T}^3)}\|\Delta u_n(t)\|_{L^2(\mathbb{T}^3)}\nonumber\\
&\leq& c\|\nabla u_n(t)\|_{L^2(\mathbb{T}^3)}^{3/2}\|\Delta u_n(t)\|_{L^2(\mathbb{T}^3)}^{3/2},\label{FSF}
\end{eqnarray}
where we used H\"older's inequality together with Agmon's  inequality \cite{Agmon} and the Poincar\'e inequality.\\ Using Young’s inequality with exponents $4$ and $4/3$ yields
$$|\langle (u_n\cdot\nabla u_n), \Delta u_n \rangle_{L^2(\mathbb{T}^3)}|\leq c\|\nabla u_n\|^6+\frac{\nu}{2}\| \Delta u_n(t)\|_{L^{2}}^2,$$
where $c$ is a positive constant that does not depend on $n$.
It turns out that
$$\frac{d}{dt}\| \nabla u_n(t)\|_{L^{2}}^2+\nu\| \Delta u_n(t)\|_{L^{2}}^2\leq c\|\nabla u_n\|^6.$$
By comparing the function $\| \nabla u_n(t)\|_{L^{2}}^2$ with the solution of the ODE:
$$\frac{dx}{dt}=cx^3,~~x(0)=\|\nabla u_0\|_{L^2}^2,$$
we infer that as long as $0\leq t<\frac{1}{2c\|\nabla u_0\|_{L^2}^4}$, the following holds 
\begin{equation}\label{AW}
\| \nabla u_n(t)\|_{L^{2}}^2\leq \underbrace{\frac{\|\nabla u_0\|_{L^2}^2}{\sqrt{1-2ct\|\nabla u_0\|_{L^2}^4}}}_{\mathcal{W}(t)}.
\end{equation}
From (\ref{AW}) and (\ref{FSF}) we now have uniform bounds on $u_n\in L^\infty([0,T_{max});H^1)$ and on $u_n\in L^2([0,T_{max});H^2)$ where $T_{max}\sim \frac{1}{\|\nabla u_0\|_{L^2}^4}$. Those uniform bounds together with $(NS_n)$ and a standard procedure allows to take the limit as $n\to\infty$ (see \cite{james} and references therein). The standard method shows that the limit $u$ is a strong solution on $[0,T_{max})$. However, what happens after time $T_{max}$ is unknown. \\
We turn now to the question of whether a local in time strong solution can be extended into a global solution. To this end, let us  make estimates directly for $u$ instead of using the Galerkin approximation. We know that $u_0$ gives rise to a strong solution that exists at least on a certain time interval $[0, T_{max} )$. On this time interval for each time $t \in(0, T_{max} )$ we take the $L^2$-inner product of $(NS)$ against $-\Delta u$, we obtain
\begin{eqnarray*}
\frac{1}{2}\frac{d}{dt}\| \nabla u(t)\|_{L^{2}}^2&+&\nu\| \Delta u(t)\|_{L^{2}}^2\\&\leq&
 |\langle (u\cdot\nabla u), \Delta u \rangle_{L^2(\mathbb{T}^3)}|\\
&\leq& \|u(t)\|_{L^\infty(\mathbb{T}^3)}\|\nabla u(t)\|_{L^2(\mathbb{T}^3)}\|\Delta u(t)\|_{L^2(\mathbb{T}^3)}.
\end{eqnarray*}
The Fourier expansion of $u(x,t)$ is given by 
$$u(x,t)=\sum_{k\in\mathbb{Z}^3}\hat{u}(k,t)e^{ikx}.$$
For a certain number $m$ (to be discussed later on), we have
\begin{eqnarray*}
\|u(t)\|_{L^\infty(\mathbb{T}^3)}&\leq&  \sum_{k\in\mathbb{Z}^3}|\hat{u}(k,t)|\\
&=&\sum_{|k|\leq m}|\hat{u}(k,t)|+\sum_{|k|> m}|\hat{u}(k,t)|.
\end{eqnarray*}
Two possible natural cases may occur. The first is when the major amount of energy at the instant $t$ is quantized in low-frequency components. This case can be represented by  the following inequality:
\begin{eqnarray}\label{Case1}
\sum_{|k|> m}|\hat{u}(k,t)|\leq \sum_{|k|\leq m}|\hat{u}(k,t)|.
\end{eqnarray}
The second case is when the major amount of energy at time $t$ is quantized in high-frequency components. That is to say: 
\begin{eqnarray}\label{Case2}
\sum_{|k|\leq m}|\hat{u}(k,t)|\leq \sum_{|k|> m}|\hat{u}(k,t)|.\end{eqnarray}
We state here the Agmon's inequality \cite{Agmon} which reads: $$\sum_{k\in\mathbb{Z}^3}|\hat{u}(k,t)|\leq c\|u(t)\|_{H^1(\mathbb{T}^3)}^{1/2}\|u(t)\|_{H^2(\mathbb{T}^3)}^{1/2}.$$ 
By Theorem \ref{regularity} we have $u\in C([0,T_{max});H^1)$ and $C((0,T_{max});H^2)$, then by a continuity argument one can always find at least a small interval of strictly positive length $[t_0,t_{1}]\subset (0,T_{max})$ on which either (\ref{Case1}) or (\ref{Case2}) occurs for any positive number $m$. Since $[t_0,t_1]$ will serve as a test interval to examine case (\ref{Case1}) and case (\ref{Case2}), one can choose without loss of generality $t_0$ very close to the instant zero. To be more precise, let $t_0$ be the instant of time immediately after $t=0$, there exists $t_1>t_0$ such that we have either (\ref{Case1}) for all $t\in [t_0,t_1]$ or (\ref{Case2}) for all $t\in [t_0,t_1]$.\\
We need also to make use of Lemma \ref{regularity2}, which combined with the fact that $\partial_t|\hat{u}(k,t)|\leq |\partial_t\hat{u}(k,t)|$ yields 
\begin{eqnarray}\label{derivative}
\partial_t\sum_{k\in\mathbb{Z}^3}|\hat{u}(k,t)|\in L^\infty(\varepsilon,T_{max};H^p),~~\forall p\in\mathbb{N}.
\end{eqnarray}
Property (\ref{derivative}) is useful because it assures the smoothness of $\sum_{k\in\mathbb{Z}^3}|\hat{u}(k,t)|$ with respect to time and hence that of $\sum_{|k|\leq m}|\hat{u}(k,t)|$ and $\sum_{|k|> m}|\hat{u}(k,t)|$. This prevents the abrupt bends of the function $t\mapsto F_m(t)=\sum_{|k|\leq m}|\hat{u}(k,t)|-\sum_{|k|> m}|\hat{u}(k,t)|$.  \\ Let us now discuss both cases on $[t_0,t_1]$. To be more precise, if (\ref{Case1}) holds true on $[t_0,t_1]$ what will happen and if (\ref{Case2}) holds true on $[t_0,t_1]$ what will happen.  \\
\textit{If condition (\ref{Case1}) holds:} \\ 
By using (\ref{Case1}), the Cauchy-Schwarz inequality and Young's inequality, we get
\begin{align*}
\frac{1}{2}\frac{d}{dt}\| \nabla u(t)\|_{L^{2}}^2&+\nu\| \Delta u(t)\|_{L^{2}}^2\\&\leq 2\sum_{|k|\leq m}|\hat{u}(k,t)|\|\nabla u(t)\|_{L^2(\mathbb{T}^3)}\|\Delta u(t)\|_{L^2(\mathbb{T}^3)}\\
&\leq 2\left(\sum_{|k|\leq m}1\right)^{1/2}\left(\sum_{|k|\leq m}|\hat{u}(k,t)|^2\right)^{1/2}\\&\times\|\nabla u(t)\|_{L^2(\mathbb{T}^3)}\|\Delta u(t)\|_{L^2(\mathbb{T}^3)}\\
&\leq 2\left(\sum_{|k|\leq m}1\right)^{1/2}\|u(t)\|_{L^2}\|\nabla u(t)\|_{L^2(\mathbb{T}^3)}\|\Delta u(t)\|_{L^2(\mathbb{T}^3)}\\
&\leq C(m)\|u(t)\|_{L^2}^2\|\nabla u(t)\|_{L^2(\mathbb{T}^3)}^2+\frac{\nu}{2}\|\Delta u(t)\|_{L^2(\mathbb{T}^3)}^2,
\end{align*}
where $C(m)=\displaystyle\frac{2\sum_{|k|\leq m}1}{\nu}$. By using the energy inequality for weak solutions (\ref{W.E}) and dropping the viscous term from both sides above, we obtain
$$\frac{d}{dt}\| \nabla u(t)\|_{L^{2}}^2\leq 2C(m)\|u_0\|_{L^2}^2\|\nabla u(t)\|_{L^2(\mathbb{T}^3)}^2.$$
The Gronwall's inequality yields 
\begin{eqnarray*}\label{easy1}
\| \nabla u(t)\|_{L^{2}}^2\leq \| \nabla u(t_0)\|_{L^{2}}^2\exp\{2C(m)\|u_0\|_{L^2}^2(t-t_0)\},~~\mbox{for all~}t\in[t_0,t_1].
\end{eqnarray*}
\textit{If condition (\ref{Case2}) holds:}\\
By using (\ref{Case2}) and the Cauchy-Schwarz inequality we infer that
\begin{eqnarray*}
\frac{1}{2}\frac{d}{dt}\| \nabla u(t)\|_{L^{2}}^2&+&\nu\| \Delta u(t)\|_{L^{2}}^2\\&\leq& 2\sum_{|k|> m}|\hat{u}(k,t)|\|\nabla u(t)\|_{L^2(\mathbb{T}^3)}\|\Delta u(t)\|_{L^2(\mathbb{T}^3)}\\
&=& 2\sum_{|k|> m}|k|^{-2}|k|^2|\hat{u}(k,t)|\|\nabla u(t)\|_{L^2(\mathbb{T}^3)}\|\Delta u(t)\|_{L^2(\mathbb{T}^3)}\\
&\leq& 2\left(\sum_{|k|> m}|k|^{-4}\right)^{1/2}\left(\sum_{|k|> m}|k|^{4}|\hat{u}(k,t)|^2\right)^{1/2}\\&\times&\|\nabla u(t)\|_{L^2(\mathbb{T}^3)}\|\Delta u(t)\|_{L^2(\mathbb{T}^3)}.
\end{eqnarray*}
We recall that $\sum_{|k|> m}|k|^{4}|\hat{u}(k,t)|^2\leq \sum_{k\in \mathbb{Z}^3}|k|^{4}|\hat{u}(k,t)|^2=\|\Delta u(t)\|_{L^2(\mathbb{T}^3)}^2$ and 
\begin{eqnarray*}
\sum_{|k|> m}|k|^{-4}&\leq& c_1\int_m^\infty\frac{\kappa^2}{\kappa^4}d\kappa\\
&\leq& c_1m^{-1}.
\end{eqnarray*}
It turns out
$$\frac{1}{2}\frac{d}{dt}\| \nabla u(t)\|_{L^{2}}^2+\{\nu-2c_1^*
{m^{-1/2}}\| \nabla u(t)\|_{L^{2}}\}\| \Delta u(t)\|_{L^{2}}^2\leq 0,$$ where $c_1^*=\sqrt{c_1}$.  Since  $\lim\limits_{m\to\infty}{m^{-1/2}}=0$, then one can choose the number $m$ such that 
 $$m>\frac{4}{c_1\nu^2}\|\nabla u(t_0)\|_{L^2(\mathbb{T}^3)}^{2}.$$
In such a way, the factor $\{\nu-2c_1^*{m^{-1/2}}\|\nabla u(t)\|_{L^2(\mathbb{T}^3)}\}$ would still positive at least over a short interval of time $[t_0,\tau_1]\subset[t_0,t_1]$.  Consequently, it turns out that
$$\frac{d}{dt}\|\nabla u(t)\|_{L^2(\mathbb{T}^3)}^2\leq 0,~~\mbox{and}~\|\nabla u(t)\|_{L^2(\mathbb{T}^3)}\leq \|\nabla u(t_0)\|_{L^2(\mathbb{T}^3)} ~\mbox{for all~}t\in [t_0,\tau_1]. $$
But as $\|\nabla u(t)\|_{L^2(\mathbb{T}^3)}$ is continuous on $[t_0,t_1]$, we obtain $$\|\nabla u(t)\|_{L^2(\mathbb{T}^3)}\leq \|\nabla u(t_0)\|_{L^2(\mathbb{T}^3)} ~\mbox{for all~}t\in [t_0,t_1].$$
Thus, the condition on $m$ has been determined successfully. In fact, by choosing such number $m$ (i.e. ${m^{-1/2}}<\frac{\nu}{2c_1^*}\|\nabla u(t_0)\|_{L^2}^{-1}$) one ensures that $\|\nabla u(t)\|_{L^2}$ is controlled on the time interval $[t_0,t_1]$ regardless of the sign of $F_m(t)=\sum_{|k|\leq m}|\hat{u}(k,t)|-\sum_{|k|> m}|\hat{u}(k,t)|$ on it. It is worth mentioning that starting from the instant $t_0$, the procedure above remains applicable  as long as $F_m(t)$ keeps its sign constant until it reverses the sign at the instant $t_1$. At the instant $t_1$, the function $F_m(t)$ should be updated by changing the condition on $m$. Precisely, we take another number $m$ such that ${m^{-1/2}}<\frac{\nu}{2c_1^*}\|\nabla u(t_1)\|_{L^2}^{-1}$.\\
Let us associate the number $m_0=\frac{8}{c_1\nu^2}\|\nabla u(t_0)\|_{L^2(\mathbb{T}^3)}^{2}$ to the interval $[t_0,t_1]$ on which $F_{m_0}(t)$ keeps its sign constant.
To conclude, we have proved that there exists $t_1>t_0$ such that we have for all $t\in[t_0,t_1]$: 
\begin{eqnarray*}
\| \nabla u(t)\|_{L^{2}}^2&\leq& \| \nabla u(t_0)\|_{L^{2}}^2\exp\{2C(m_0)\|u_0\|_{L^2}^2(t-t_0)\}
\\
&\leq& \| \nabla u(t_0)\|_{L^{2}}^2\exp\left\{\frac{2048\times c_2}{c_1^3\nu^7}\|\nabla u(t_0)\|_{L^2(\mathbb{T}^3)}^{6}\|u_0\|_{L^2}^2(t-t_0)\right\},
\end{eqnarray*}
if condition (\ref{Case1}) holds true on $[t_0,t_1]$. We point out that we used $$\sum_{|k|\leq m}1\leq c_2\int_0^mk^2d\kappa=c_2m^3=c_2\frac{8^3}{c_1^3\nu^6}\|\nabla u(t_0)\|_{L^2(\mathbb{T}^3)}^{6}.$$
Or
\begin{eqnarray*}
\| \nabla u(t)\|_{L^{2}}^2\leq \| \nabla u(t_0)\|_{L^{2}}^2,
\end{eqnarray*}
if condition (\ref{Case2}) holds true on $[t_0,t_1]$.
The statement above can be summarized by making use of the fact that $$\| \nabla u(t_0)\|_{L^{2}}^2\leq \| \nabla u(t_0)\|_{L^{2}}^2\exp\{\frac{2048\times c_2}{c_1^3\nu^7}\|\nabla u(t_0)\|_{L^2(\mathbb{T}^3)}^{6}\|u_0\|_{L^2}^2(t-t_0)\}.$$ Thus, there exists $t_1>t_0$ such that we have for all $t\in[t_0,t_1]:$  $$\| \nabla u(t)\|_{L^{2}}^2\leq \| \nabla u(t_0)\|_{L^{2}}^2\exp\left\{\frac{2048\times c_2}{c_1^3\nu^7}\|\nabla u(t_0)\|_{L^2(\mathbb{T}^3)}^{6}\|u_0\|_{L^2}^2(t-t_0)\right\}.$$
Continuing in this vein, in the next interval we know already that $m$ must be as large as ${m^{-1/2}}<\frac{\nu}{2c_1}\|\nabla u(t_1)\|_{L^2}^{-1}$ which guarantees by continuity that in case (\ref{Case2}) the function $\| \nabla u(t)\|_{L^{2}}$ is non-increasing on this interval. There exists then $t_2>t_1$ such that for all $t\in[t_1,t_2]$:
\begin{eqnarray*}
\| \nabla u(t)\|_{L^{2}}^2\leq \| \nabla u(t_1)\|_{L^{2}}^2\exp\left\{\frac{2048\times c_2}{c_1^3\nu^7}\|\nabla u(t_1)\|_{L^2(\mathbb{T}^3)}^{6}\|u_0\|_{L^2}^2(t-t_1)\right\},
\end{eqnarray*}
if condition (\ref{Case1}) holds true on $[t_1,t_2]$.
Or
\begin{eqnarray*}
\| \nabla u(t)\|_{L^{2}}^2\leq \| \nabla u(t_1)\|_{L^{2}}^2,
\end{eqnarray*}
if condition (\ref{Case2}) holds true on $[t_1,t_2]$.
Repeating this process as many times as needed to obtain $[t_0,T_{max})=\cup_{j=0}^{N-1}[t_j,t_{j+1}]\cup [t_{N},T_{max})$ (where $\epsilon$ is an arbitrary small constant and $[t_j,t_{j+1}]$ are successive intervals), such that for all $t\in [t_j,t_{j+1}]$ we have either 
\begin{eqnarray}\label{F.S.F1}
\| \nabla u(t)\|_{L^{2}}^2\leq \| \nabla u(t_j)\|_{L^{2}}^2\exp\left\{\frac{2048\times c_2}{c_1^3\nu^7}\|\nabla u(t_j)\|_{L^2(\mathbb{T}^3)}^{6}\|u_0\|_{L^2}^2(t-t_j)\right\}
\end{eqnarray}
or
\begin{eqnarray}\label{F.S.F2}
\| \nabla u(t)\|_{L^{2}}^2\leq \| \nabla u(t_j)\|_{L^{2}}^2.
\end{eqnarray}
This process would certainly control the norm $\| \nabla u(t)\|_{L^{2}}^2$ and rules out the blowup of $u$ in $H^1(\mathbb{T}^3)$ as $t$ approaches $T_{max}$. In fact, on the interval $[t_N,T_{max})$ either $F_{m_N}(t)\geq 0$ holds true for all $t\in [t_N,T_{max})$ or $F_{m_N}(t)\leq 0$ holds true for all $t\in [t_N,T_{max})$ where $m_N=\frac{8}{c_1\nu^2}\|\nabla u(t_N)\|_{L^2(\mathbb{T}^3)}^{2}$ and hence \begin{eqnarray*}\lim\limits_{t\to T_{max}}\|\nabla u(t)\|_{L^2(\mathbb{T}^3)}^2\leq \ \ \ \ \ \ \ \ \ \ \ \ \ \ \ \ \ \ \ \ \ \ \ \ \ \\ \|\nabla u(t_N)\|_{L^2(\mathbb{T}^3)}^2\exp\left\{\frac{2048\times c_2}{c_1^3\nu^7}\|\nabla u(t_N)\|_{L^2(\mathbb{T}^3)}^{6}\|u_0\|_{L^2}^2(T_{max}-t_N)\right\}.\end{eqnarray*}
As $u(t_N)\in H^1(\mathbb{T}^3)$, then the upper bound $$ \|\nabla u(t_N)\|_{L^2(\mathbb{T}^3)}^2\exp\left\{\frac{2048\times c_2}{c_1^3\nu^7}\|\nabla u(t_N)\|_{L^2(\mathbb{T}^3)}^{6}\|u_0\|_{L^2}^2(T_{max}-t_N)\right\}$$ is finite.
Therefore, the solution $u$ can be extended into a global in time strong solution.
\section{Discussion}
An interesting observation is the following: it suffice that (\ref{Case2}) occurs only once on an interval of strictly positive length $[t_0,t_1]\subset (0,T_{max})$ for $m$ such that $m^{-1/2}<\frac{\nu}{2c_1^*}\|\nabla u(t_0)\|_{L^2}^{-1}$ to extend the solution onto a larger interval of time $[0,T_{max}+t_1-t_0)$.
In fact, estimate (\ref{AW}) tells us that an initial data as large as $\|\nabla u_0\|_{L^2}$ gives rise to a solution $u$ that would remain bounded on an interval $[0,T_{max})$ of length $T_{max}-0$. Let $t_0\in (0,T_{max})$, inequality (\ref{AW}) also tells us that an initial data as large as $\frac{\|\nabla u_0\|_{L^2}}{\{1-2ct_0\|\nabla u_0\|_{L^2}^4\}^{1/4}}$ gives rise to a solution $u$ that remains bounded on an interval of length $T_{max}-t_0$. 
According to the proof given in the previous section, if condition (\ref{Case2}) holds true on $[t_0,t_1]$, it turns out
$$\frac{d}{dt}\|\nabla u(t)\|_{L^2(\mathbb{T}^3)}^2\leq 0,~~\mbox{and}~\|\nabla u(t)\|_{L^2(\mathbb{T}^3)}\leq \|\nabla u(t_0)\|_{L^2(\mathbb{T}^3)} ~\mbox{for all~}t\in [t_0,t_1].$$
But as we have: $$\|\nabla u(t_1)\|_{L^2}^2\leq \|\nabla u(t_0)\|_{L^2}^2\leq \frac{\|\nabla u_0\|_{L^2}^2}{\sqrt{1-2ct_0\|\nabla u_0\|_{L^2}^4}},$$ therefore by starting from $t_1$ the solution $u$ would now still bounded on an interval of time of length $T_{max}-t_0$. In other words, the solution $u$ is extended to the interval $[0,T_{max}+t_1-t_0)$ which means that $u(T_{max})\in H^1(\mathbb{T}^3)$.

At this point, one may ask the question under which condition the norm $\|\nabla u(t)\|_{L^2}$ keeps decreasing for all positive time $t\in \mathbb{R}_+.$ In fact, this is possible when the distribution of energy in the initial data is extremely unbalanced (i.e. $\sum_{|k|\leq m}|\hat{u}(k,0)|<<\sum_{|k|> m}|\hat{u}(k,0)|$). In that case, by smoothness of the function  $\sum_{k\in\mathbb{Z^3}}|\hat{u}(k,t)|$ with respect to time, condition (\ref{Case2}) keeps for a long interval of time until potentially $\|\nabla u(t)\|_{L^2}$ satisfies the smallness condition of \cite{kato}. \\ Another aspect to discuss here is the motivation behind choosing the instant $t_0$ very close to zero. In fact, by doing so one can ensure via (\ref{AW}) the closeness of $\|\nabla u(t_0)\|_{L^2}$ to $\|\nabla u_0\|_{L^2}$ while holding the necessary regularity ($u(t_0)\in H^p(\mathbb{T}^3)$ for all $p\in\mathbb{N}$). This also guarantees the minimum worsening to $\|\nabla u(t_1)\|_{L^2}$ in case (\ref{Case1}). However, it is needless to say that this was optional and that one can choose any instant $t_0\in[0,T_{max})$ as initial time.

\section{Conclusion}

We have already proved that the local in time strong solution to ($NS$) can be extended to become global in time strong solution. This was done via making estimates to $u$ in $H^1$ on a series of time intervals requiring that the function 
$t\mapsto F_m(t)=\sum_{|k|\leq m}|\hat{u}(k,t)|-\sum_{|k|> m}|\hat{u}(k,t)|$ keeps its sign constant (either positive or negative) on each of them. It is worth noting that when the major amount of energy is located in high-frequency components (i.e. $F_m(t)\leq 0$), the norm $\|\nabla u(t)\|_{L^2}$ decreases with time. This is in fact consistent with the phenomenology of the turbulent cascade which states that energy is dissipated at the small scales (i.e. higher frequencies). \\

\end{document}